\documentclass[letterpaper, 10 pt, conference]{ieeeconf}
\IEEEoverridecommandlockouts  
\usepackage{lineno,setspace, nicefrac, subcaption, caption}
\usepackage{mathrsfs}
\usepackage{tikz}
\usetikzlibrary{automata,arrows,positioning,calc}
\usepackage[mathscr]{eucal}
\usepackage{float}
\usepackage{cite}
\usepackage[hidelinks]{hyperref}
\DeclareGraphicsExtensions{.png,.jpeg,.jpg,.pdf}
\graphicspath{{./Figures/}}

\newtheorem{theorem}{Theorem}[section]
\newtheorem{definition}[theorem]{Definition} 

\newtheorem{corollary}[theorem]{Corollary}
\newtheorem{lemma}[theorem]{Lemma}

\newtheorem{remark}[theorem]{Remark}
\newtheorem{notation}[theorem]{Notation}

\newcommand{\lb}{\left\{}
\newcommand{\rb}{\right\}}
\newcommand{\set}[1]{\left\{ #1 \right\}}

\newcommand{\Def}{\overset{\textbf{def}}{=}}
\newcommand{\RR}{\mathbb{R}}

\newcommand{\NN}{\mathbb{N}}

\newcommand{\eps}{\varepsilon}

\newcommand{\Bv}{\mathbf v}
\newcommand{\Bx}{\mathbf x}
\newcommand{\abs}[1]{\left \lvert #1 \right \rvert}
\newcommand{\norm}[1]{\left\| #1 \right\|}

\newcommand{\sT}{\mathsf T}

\usepackage{amsopn}

\DeclareMathOperator{\loc}{loc}

\newcommand{\sa}{\mathsf a}

\newcommand{\udelta}{\underline{\delta}}
\newcommand{\domain}{\mathcal S}
\newcommand{\By}{\mathbf y}

\newcommand{\Bw}{\mathbf{w}}
\newcommand{\odelta}{\bar \delta}

\usepackage{mathptmx} 
\usepackage{times} 
\usepackage{amsmath} 
\usepackage{amssymb}  

\usepackage{color}
\definecolor{darkgreen}{rgb}{0,0.5,0}
\usepackage[colorinlistoftodos]{todonotes}
\usepackage{soul}
\setulcolor{red} 
\setstcolor{red} 
\sethlcolor{yellow} 

\newcount\Comments  
\usepackage{color}
\definecolor{darkgreen}{rgb}{0,0.5,0}
\definecolor{purple}{rgb}{1,0,1}
\newcommand{\kibitz}[2]{\ifnum\Comments=0\textcolor{#1}{#2}\fi}

\usepackage{dirtytalk}

\title{\LARGE \bf {On the Analytical Properties of a Nonlinear Microscopic Dynamical Model for Connected and Automated Vehicles}}
\author{Hossein Nick Zinat Matin$^{1}$, Yuneil Yeo$^{1}$, Xiaoqian Gong$^{2}$ and {Maria Laura} {Delle Monache}$^{1}$ 
\thanks{*Y.Yeo is partially supported by the Dwight David Eisenhower Transportation Fellowship Program Under Grant No. 693JJ32445085.}
\thanks{$^{1}$ H.N.Z.Matin, Y.Yeo and M.L.Delle Monache are with the Department of Civil and Environmental Engineering,
        University of California, Berkeley,
        {\tt\small h-matin@berkeley.edu, yuneily@berkeley.edu, mldellemonache@berkeley.edu.}}    
\thanks{$^{2}$ X.Gong is with the Department of Mathematics and Statistics, Amherst College, Amherst, Massachusetts, 
{\tt\small xgong@amhert.edu}.}
}

\date{\today.}

\begin{document}

\maketitle 
\thispagestyle{empty}
\pagestyle{empty}

\begin{abstract}
In this paper, we propose an integrated dynamical model of Connected and Automated Vehicles (CAVs) which incorporates CAV technologies and a microscopic car-following model to improve safety, efficiency and convenience. We rigorously investigate the analytical properties such as well-posedness, maximum principle, perturbation and stability of the proposed model in some proper functional spaces. Furthermore, we prove that the model is collision free and we derive and explicit lower bound on the distance as a safety measure.
\end{abstract}
\section{Introduction and Related Works}
Traffic flow models can be categorized into microscopic, mesoscopic, and macroscopic models depending on the scale at which traffic is represented \cite{treiber2013traffic}. In this article, we focus on microscopic models that describe the interaction between the individual vehicles as well as with the mixed autonomy condition. 
Over the past decades, various works contributed to microscopic models \cite{pipes1953operational,newell1961nonlinear, gazis1961nonlinear,treiber2017intelligent}, including 
optimal velocity model (OVM) \cite{bando1998analysis, wang2023car} which cannot prevent collision, the intelligent driver model (IDM) \cite{treiber2017intelligent, zha2023survey} which are proven to be unsuccessful to capture the characteristics of CAVs, \cite{milanes2013cooperative}, and desired measure models \cite{kerner2015failure} which ignore the communication capability of the CAVs \cite{jafaripournimchahi2022integrated}.

Car-following based dynamical models \cite{chandler1958traffic}, e.g. Optimal Velocity Follow-the-Leader (OVFL), employ the interaction between the singularity term (the inverse of the distance $(\Delta x)^{-1}$), and relative velocity, $\Delta v$, to ensure a collision-free and relatively stable behavior of solution \cite{nick2022near,dellemonache2019pardalos, gong2023well}. These models are considered as a suitable framework for CAV dynamics, however, they ignore the \textit{desired velocity} (applied by drivers or by centralized control) which makes them less efficient in \textit{mixed autonomy condition}, \cite{matin2023existence}. 

With the emergence of autonomous driving technologies, Advanced Driving-Assistant Systems (ADAS) such as Adaptive Cruise Control (ACC), Cooperative Adaptive Cruise Control (CACC), and self-driving systems, several works expanded microscopic models to study the behavior of autonomous vehicles (AVs) with the goal of improving comfort and safety for drivers, \cite{van2006impact, milanes2014modeling, shladover2012impacts}. While the existing models take advantage of communication capabilities of CAVs, they do not provide a \textit{collision free framework}.

In this work, we propose a novel integrated nonlinear dynamical model of CAVs that provides a framework for both mixed autonomy condition (where AVs and human-driven vehicles coexist) and CAV platooning (where AVs dynamically interact).
The proposed model extends the existing literature in several directions.

Firstly, in contrast to the OVFL-based models \cite{nick2022near, matin2023near}, by integrating a generic control term, our model takes the desired velocity into account. This creates a more flexible framework from \textit{communication stand point} for a broad range of control applications; e.g. the designing control in micro-macro presentation of \textit{mixed autonomy} \cite{matin2023existence}, \cite{wang2024hierarchical}, or acting as an adaptive cruise control (cf. \cite{van2006impact}) to regularize the velocity profile.

Next, \textit{from safety point of view}, unlike the CACC and ACC models, \cite{van2006impact}, our proposed model will be rigorously proven to be \textit{collision free} and will \textit{not experience negative velocity}. In addition, while most of the analytical results in the literature are from control, stability and simulation standing points, \cite{davis2003modifications, stern2018dissipation}, we employ a rigorous framework which allows us to carefully study the effect of the singularity in near-collision regions.
Such a careful investigation contributes to (i) designing efficient controls for a \textit{safe and comfortable transition between the states}, and (ii) analyzing the behavior of the system in real scenarios in which any physical system could be perturbed inevitably into different states due to various perturbation forces. Using a novel approach, we show that collision is precluded in the proposed dynamical model which proves the efficiency of the dynamical model from the safety point of view. More importantly, we derive an explicit lower bound on the distance as a safety measure.

Finally, in this work, we prove several key properties of the proposed model which are crucial in understanding the behavior of the solution. More precisely, we propose a novel and rigorous analysis of well-posedness and construct a unique solution. The maximum principle for velocity is proven through the interaction of the singularity and control terms. A careful stability study of the system by rigorously analyzing the perturbed dynamics has been investigated. Such analyses provide a helpful framework for further investigation of this type of dynamical systems.

The organization of this paper is as follows. We introduce the mathematical model of the problem in Section \ref{S:Introduction}. In Section \ref{S:well-posedness}, the analytical properties of the model are elaborated. A small perturbation analysis of the system and perturbed system's convergence to the original dynamics is discussed in Section \ref{S:perturbation}. Some preliminary estimates on the solutions along the trajectories are derived in Section \ref{S:Estimates}.
\section{Mathematical Model}\label{S:Introduction}
 Let $T>0$ to be any fixed time horizon. For {$N \in \NN$, and $N \ge 2$}, we {propose} the dynamical model of the form:
\begin{equation}\label{E:general}
    \begin{cases}
        \dot x_1 = v_1\\
        \dot v_1 = \sa_1\\
        \dot x_{n} = v_n\\
        \dot v_n = \min\set{k_v \frac{( v_{n -1} -  v_n)}{(x_{n-1} - x_n)^2} + k_d (x_{n -1} - x_n - \tau_s v_n), k (u_n - v_n)}\\
        (x_n(0), v_n(0)) = (x_{n,\circ}, v_{n,\circ}),
    \end{cases}
\end{equation}
for $n \in \{2, \dots, N\}$. The acceleration term $\sa_1$ of the first vehicle is assumed to belong to $C\cap L^{\infty}([0,T]; \mathbb{R})$, $x_n$ and $v_n$ are the position and velocity of the $n$-th vehicle, respectively.
Here, $k_v, k_d, k$ and $\tau_s >0$ are constant weights that can be calibrated to reflect the desired intensity of each term. Functions $u_n \in C([0, T]; [\underline u, \bar u])$ are the controls (can be interpreted as desired velocity applied by drivers or centralized control) where $0 < \underline u < \bar u < \bar v$. The constant $\bar v$ is the maximal velocity. The interval $[\underline u, \bar u]$ is the range of the controls and are used for analytical purposes.
In the case of our proposed dynamics, it is sufficient to analyze the system for $N =2$, i.e., two consecutive vehicles, and the generalization of the results to $N>2$ will be immediate (see Remark \ref{R:generalization}). Therefore, we consider the reduced dynamics in the form of: 
\begin{equation}\label{E:dynamics}
    \begin{cases}
        \dot x_\ell = v_{\ell} \\
        \dot v_\ell = \sa_\ell \\
        \dot x = v\\
        \dot v = \mathcal A(x_\ell, x, v_\ell, v; u)\\
        (x_\ell(0), x(0), v_\ell(0),v(0)) = (x_{\ell, \circ}, x_\circ, v_{\ell, \circ}, v_\circ),
        \end{cases}
\end{equation}
and we consider
\begin{equation}\label{E:leading_accel}
   0\leq v_{\ell,\circ}+ \int_0^t \sa_{\ell}(s) ds \le \bar v < \infty,  \forall t \in [0, T] .
\end{equation}
In particular, Eq. \eqref{E:leading_accel} implies that the leading vehicle is moving within the speed limit and with bounded acceleration. 
Moreover, the acceleration map $(x_\ell, x, v_\ell, v) \in \RR^2 \times \RR_+^2 \mapsto \mathcal A(x_\ell, x, v_\ell, v; u) \in \RR$ is defined by
\begin{equation}\label{E:accel}\begin{split}
  &\mathcal A(x_\ell, x, v_\ell, v; u) \Def\\
   & \qquad \min\set{k_v \frac{( v_\ell -  v)}{(x_\ell- x)^2} + k_d (x_\ell- x - \tau_s v), k (u - v)}.
\end{split}\end{equation}

\begin{figure}[h!] 
  \begin{subfigure}[b]{0.49\linewidth}
    \centering
    \includegraphics[width=0.8\linewidth]{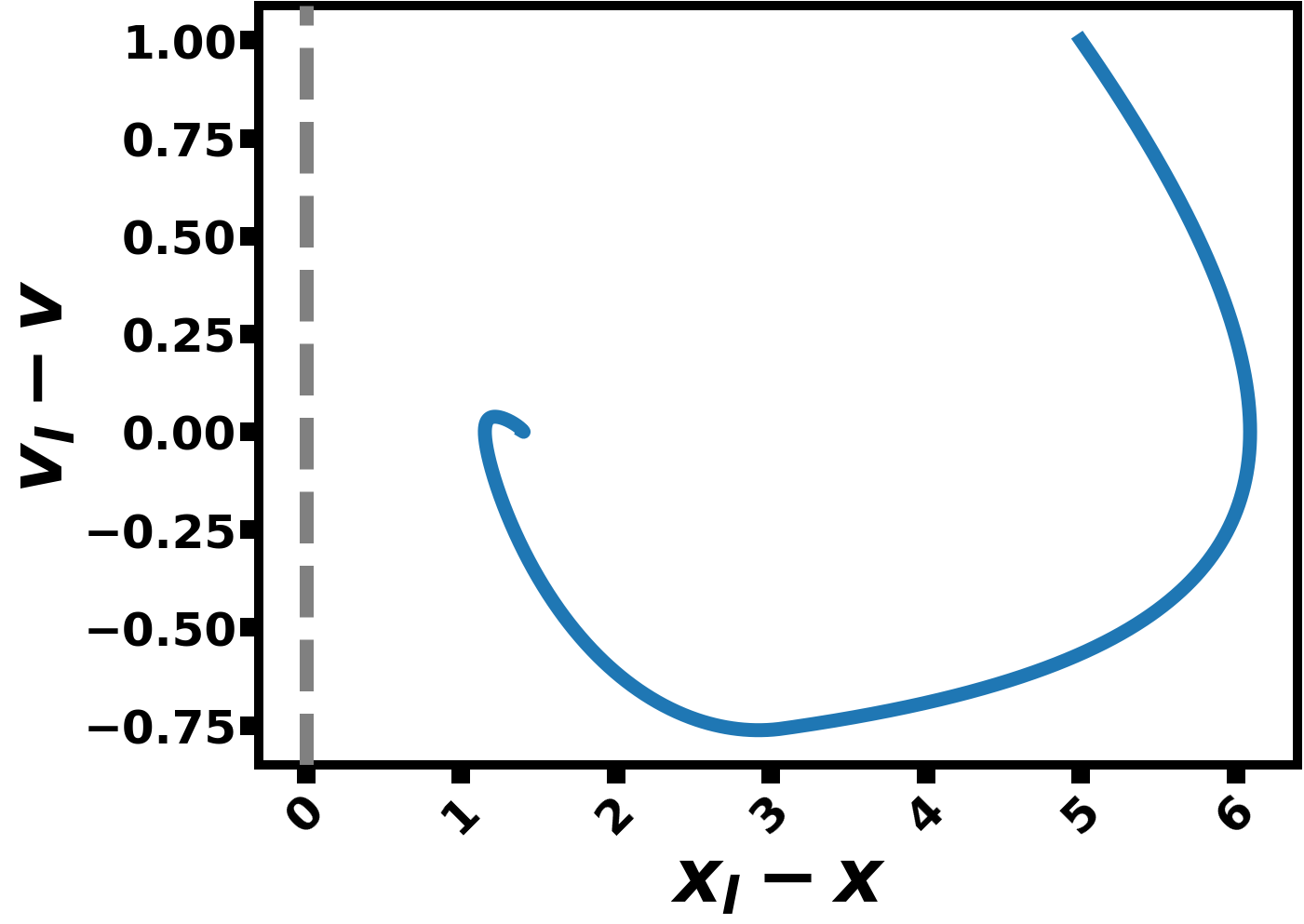} 
    \label{fig: right_top_constant_velocity} 
  \end{subfigure} 
  \begin{subfigure}[b]{0.45\linewidth}
    \centering
    \includegraphics[width=0.8\linewidth]{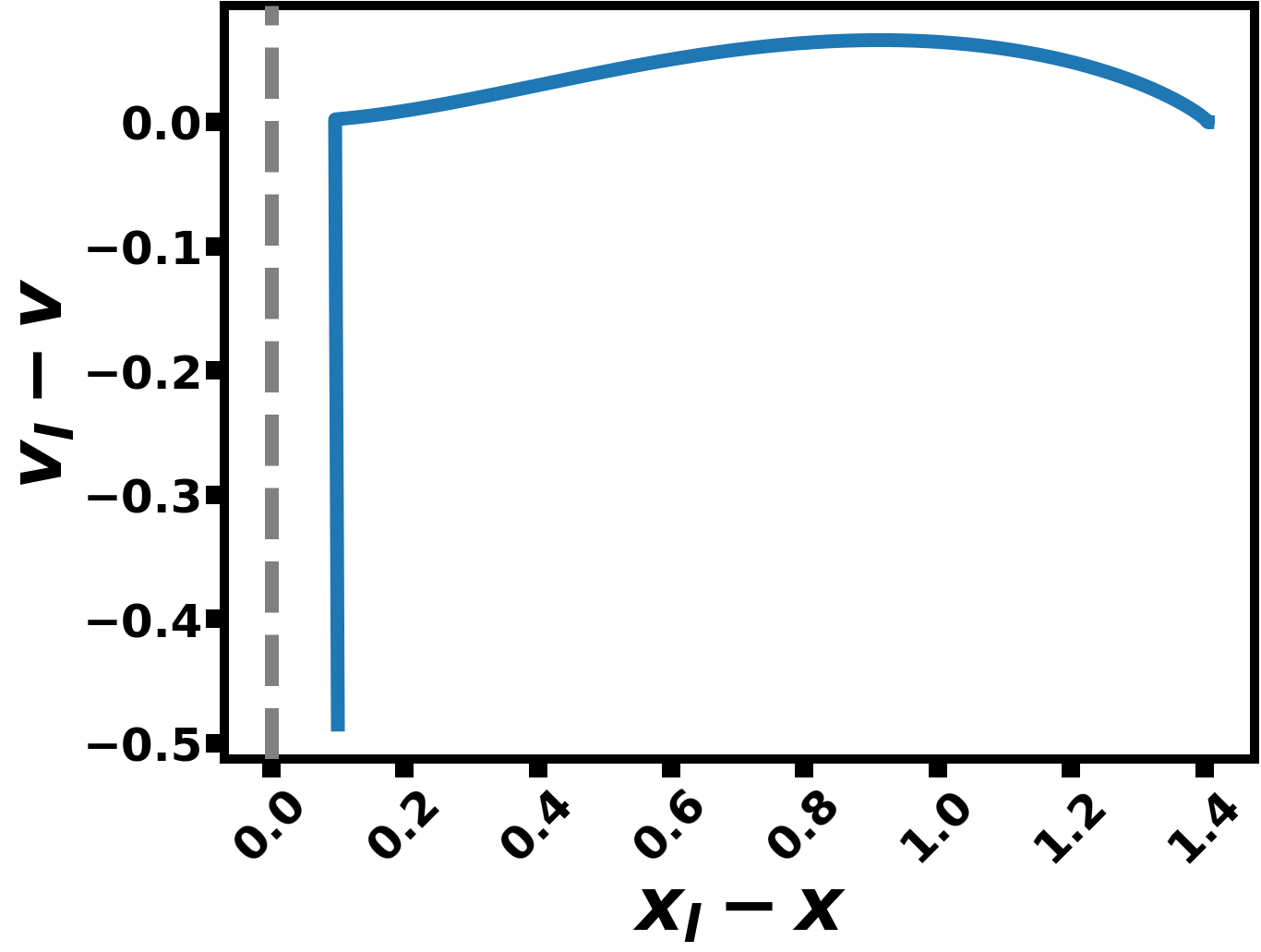} 
  \end{subfigure}
  \vspace{-5pt}
  \caption{\scriptsize{Trajectory of the Proposed Dynamics. Parameters: $k_v$ =1, $k_d$ = 0.2, $k$ = 0.3, $\tau$ = 1.4, $u$ = 1.9, $T$ = 100, $\sa_\ell$ = 0. Initial Conditions: (Left) $(x_{\ell, \circ}, x_\circ, v_{\ell, \circ}, v_\circ)$ = (5,0,1,0);  (Right) $(x_{\ell, \circ}, x_\circ, v_{\ell, \circ}, v_\circ)$ = (0.1,0,1,1.485). The gray line shows the collision.}}
 
  \label{fig: proposed} 
\end{figure}

Figure \ref{fig: proposed} shows the trajectory of the proposed dynamics with two different initial data.
\begin{remark} Before we move to the analytical results, let us briefly elaborate on the motivation of the proposed model by comparing \eqref{E:accel} with the corresponding terms in CACC and OVFL dynamics. 
We recall that in the CACC model, \cite{van2006impact}: 
\begin{equation}\label{E:acc_old}
    \dot v = \min\Big\{k_a \sa_\ell + k_v(v_\ell-v)+k_d (x_\ell - x - \Gamma(v)), k(u - v)\Big\},
\end{equation}
where, $  \Gamma(v) \Def \max \set{2, \left(\frac{1}{d} - \frac{1}{d_\ell} \right) v^2, \tau_s v}$, $d$ and $d_\ell$ are constants explaining the acceleration/deceleration capabilities of the corresponding vehicles. For the OVFL dynamics
\begin{equation*}
    \dot v = k_v \frac{v_\ell - v}{(x_\ell - x)^2} + k_d \lb V(x_\ell - x) - v \rb,
\end{equation*}
where $V(x) = \tanh(x - 2) + \tanh(2)$, $x \in \RR$ is known as optimal velocity function, \cite{nick2022near}.
Figure \ref{fig:comparison} compares the proposed, OVFL and CACC for two sets of initial data.
\begin{figure}[!ht] 
  \begin{subfigure}[b]{0.49\linewidth}
    \centering
    \includegraphics[width=0.85\linewidth]{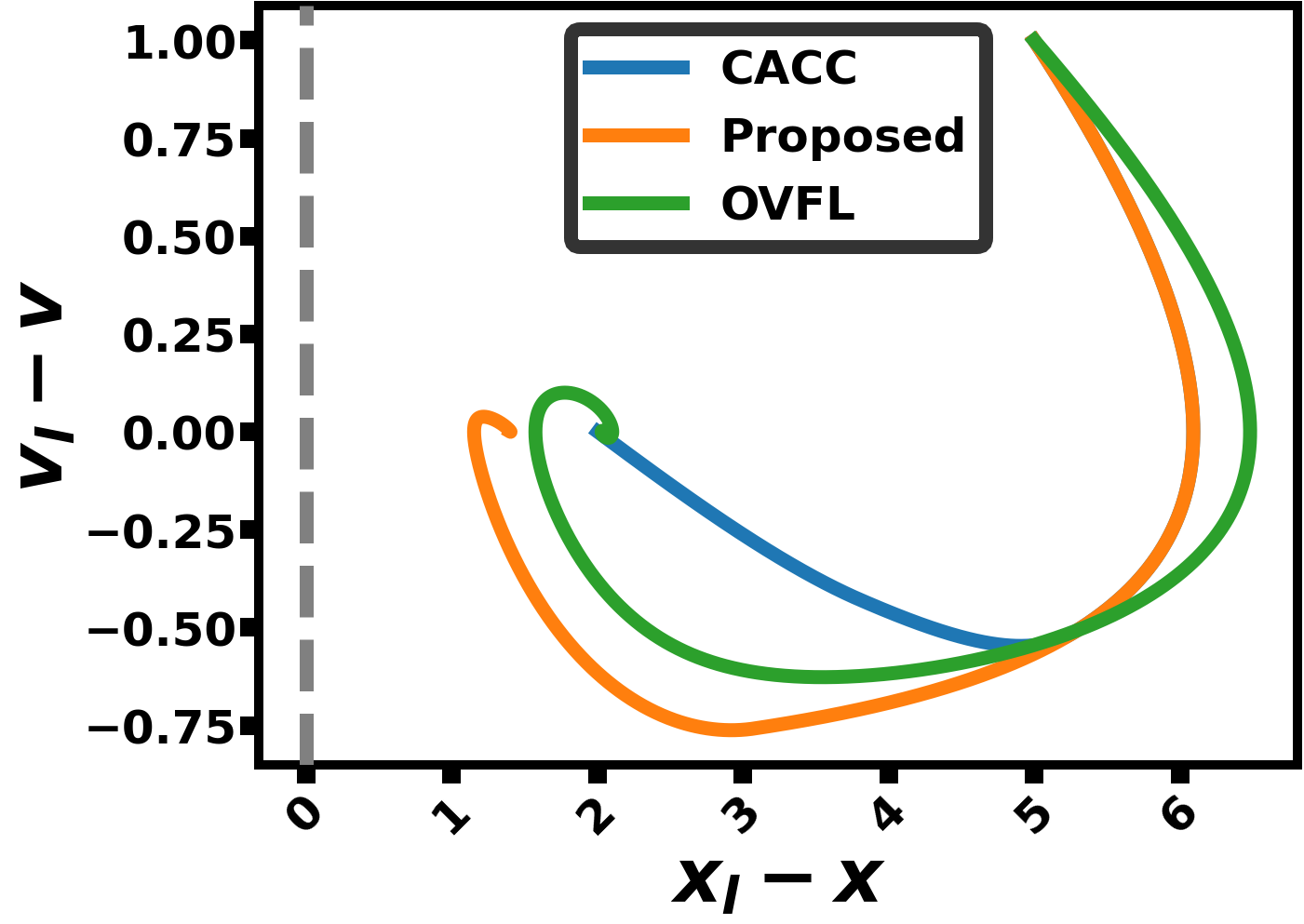} 
    \label{fig: right_top_constant_velocity} 
  \end{subfigure} 
  \begin{subfigure}[b]{0.45\linewidth}
    \centering
    \includegraphics[width=0.85\linewidth]{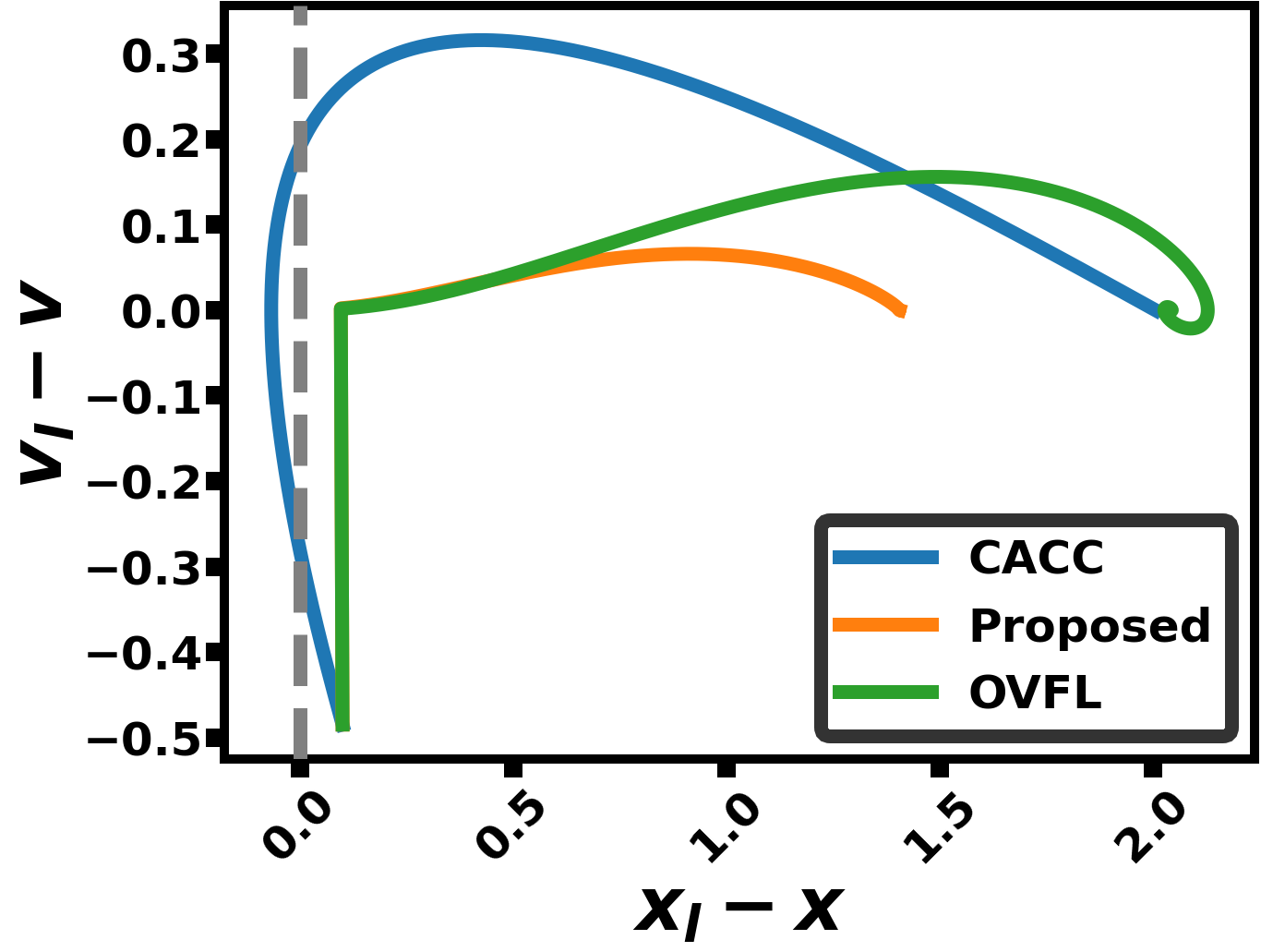} 
  \end{subfigure}
  \vspace{-5pt}
  \caption{\scriptsize{Trajectories of proposed, OVFL, and CACC. Parameters and Initial Conditions same as Fig. \ref{fig: proposed}. The gray line shows the collision. Observations from numerical experiments show that depending on the initial condition, the proposed model enjoys faster convergence to the equilibrium point in comparison with OVFL and CACC.}  
  }
  \label{fig:comparison} 
  \vspace{-10pt}
\end{figure}

The right plot of Figure \ref{fig:comparison} shows that CACC model \textit{cannot prevent collision}. More importantly, CACC can \textit{create negative velocity} at extreme cases. Such a critical issue in this model can be traced back to $\Gamma(v)$. More precisely, considering the case where the distance $x_\ell - x$ as well as the velocity function $v$ are sufficiently small, the acceleration term \eqref{E:acc_old} can encourage deceleration due to the fact that $\Gamma(v) \ge 2$, which results in negative velocity. In the proposed model, we modify the dynamics of acceleration to prevent collision and negative velocity. \hfill $\diamond$
\end{remark}
\section{Analytical Results}\label{S:well-posedness} 
 In this section, the analytical properties of the dynamical model \eqref{E:dynamics} will be elaborated. Let us start by introducing some notations.
\begin{notation}\label{def:sobolev}
    Let $\mathcal U \subset \RR_+$ and $\mathcal W \subset \RR$, $u \in L^1_{\loc}(\mathcal U; \mathcal W)$, $1 \le p < \infty$ and $1 \le q \le \infty$. Then, $u \in W^{p, q}(\mathcal U; \mathcal W)$ if $\mathcal D^\alpha u \in L^q(\mathcal U; \mathcal W)$ for $\alpha \in \set{0, \cdots, p}$, where $\mathcal D^k u$ is the $k$-th weak derivative of $u$; see \cite{adams2003sobolev} for more detail on Sobolev spaces.

    In this paper, we consider $\ell^2$-norm for vectors. For a set $\Omega \subset \RR^n$, $n \in \NN$, $C(\Omega)$ and $C^1(\Omega)$ denote the space of continuous and continuously differentiable functions, respectively. By $AC(\mathcal U; \mathcal W)$ we denote the class of absolutely continuous functions on $\mathcal U$ with values in $\mathcal W$. The boundary of a set $\Omega$ is denoted by $\partial \Omega$.  $\hfill \diamond$
\end{notation}
\begin{definition}\label{def:solution}
    Given $T>0$, $\bar v \in \RR_+$, and the initial condition $(x_{\ell, \circ}, x_\circ, v_{\ell, \circ}, v_\circ) \in \mathcal S$, with the state-space defined by 
    \small{
    \begin{equation}
        \label{E:state_space}
        \mathcal S \Def \set{(x_1, x_2, v_1, v_2) \in \RR^2 \times \RR_+^2: x_1 - x_2 > 0, v_1, v_2 \in [0, \bar v]}
    \end{equation}
    }
    \normalsize
    the $(x_\ell, x, v_\ell, v) \in W^{1, \infty}((0, T); \mathcal S)$ is the solution of \eqref{E:dynamics}, if for $t \in [0, T]$ satisfies Eqns. in \eqref{E:dynamics} in the integral form.
    Note that the Definition \ref{def:solution} implies $(x_\ell, x) \in W^{2, \infty}((0,T); \RR^2)$. 
\end{definition}
\begin{remark}[\textbf{Justification of the solution space}] The solution to \eqref{E:dynamics} in the classical sense only locally exists, and the extension of the solution to a larger domain will lose some regularity. In particular, as reflected in the Definition \ref{def:solution}, for the vector-valued solution $t \mapsto (x_\ell, x, v_\ell, v)(t)$ of \eqref{E:dynamics}, we expect to show $(x_\ell, x) \in C^1((0, T); \RR^2)$, while $(v_\ell, v) \in AC([0, T]; [0, \bar v]^2)$ and hence a.e. differentiable. In other words, in this paper, the solution of \eqref{E:dynamics} is understood in the generalized (integral) sense.\\
In addition, the function spaces provide a proper framework to understand the behavior of the solution in more general state-spaces. The choice of $W^{1, \infty}(0, T)$ is justified by the application, in particular the boundedness of the solutions.\hfill $\diamond$
\end{remark}
\begin{theorem}[\textbf{Well-posedness}]\label{T:main} 
We consider a fixed $T>0$. The dynamical model \eqref{E:dynamics} has a solution in the sense of the Definition \ref{def:solution}. In addition,
    \begin{equation}\label{E:non_collision}
        \norm{x_\ell - x}_{L^\infty(0 ,T)} \ge \frac{k_v}{v_\circ+Hk_d+\frac{k_v}{h_{\circ}}}>0,
    \end{equation}
for some constant $H = H(T)>0$ where, $h_\circ = x_{\ell, \circ} - x_\circ$. 
\end{theorem}
\begin{proof}
For coherency of the presentation, we need to prove several results along the way. The main idea is first to show that the solution exists in a specific domain uniquely. Then, using proper analytical tools we gradually extend the solution to larger domains.
Let $t = t_\circ$ be the initial time. Fix a constant $\udelta \in (0, \min \set{1, x_{\ell, \circ} - x_\circ})$, such that $(x_{\ell, \circ}, x_\circ, v_{\ell, \circ}, v_\circ) \in \Xi^{(\udelta)}$, where we define the compact set
\begin{multline}\label{E:compactSet}
    \Xi^{(\udelta)} \Def \Big\{(x_1, x_2, v_1, v_2) \in \RR^2 \times \RR_+^2: \\
    x_1 - x_2 \in [\udelta, \nicefrac{1}{\udelta}], 
    x_2 \in [-\nicefrac{1}{\udelta}, \nicefrac{1}{\udelta}], 
    v_1, v_2 \in [0, \bar v]
    \Big\}, 
\end{multline}
i.e., the initial data is located in a compact domain. In addition, for analytical purposes, let us consider the vector presentation of \eqref{E:dynamics}: 
    \begin{equation}\label{E:vec_dynamics}
        \begin{pmatrix}
            \dot \Bx \\
            \dot \Bv
        \end{pmatrix}
        = \mathcal B(t, \Bx, \Bv; u) \Def \mathcal B(t,x_\ell, x, v_\ell, v; u),
    \end{equation}
where the vector-valued function $(\Bx, \Bv) \in \mathcal S \mapsto \mathcal B(t, \Bx, \Bv; u) \in \RR^4$ is defined by the right hand side of \eqref{E:dynamics}. 
The vector-valued function $\mathcal B$ is Lipschitz continuous on $\Xi^{(\udelta)}$ since given any two functions $f_1$ and $f_2$, $  \min\set{f_1, f_2} = \frac12 \left(f_1 + f_2 - \abs{f_1 - f_2} \right)$. Therefore, if $f_1$ and $f_2$ are Lipschitz continuous, then $\abs{f_1 - f_2}$ and hence $\min\set{f_1, f_2}$ will be so. For any two points $(\Bx, \Bv), (\By, \Bw) \in \Xi^{(\udelta)}$, by properties of $\sa_\ell$ (see \eqref{E:leading_accel}) we can show that $t \mapsto \mathcal B(t, \Bx, \Bv;u)$ is continuous and
    \small{\begin{equation*}
        \norm{\mathcal B(t,\Bx, \Bv; u) - \mathcal B(t,\By, \Bw; u)}_{\ell^2(\RR^4)} \le 
        K_{\Xi^{(\udelta)}} \norm{(\Bx, \Bv) - (\By, \Bw)}_{\ell^2(\RR^4)}
    \end{equation*}
    }
    \normalsize
for some (Lipschitz) constant $K_{\Xi^{(\udelta)}}$ and where the inequality follows by compactness of $\Xi^{(\udelta)}$.\\
Let $(t_\circ, \Bx_\circ, \Bv_\circ) \in [0, T) \times int(\Xi^{(\udelta)})$, interior of set $\Xi^{(\udelta)}$, be the initial condition. By Picard-Lindel\"of theorem, the solution $(\Bx^{(\udelta)}, \Bv^{(\udelta)}) \in C \left([t_\circ, t_\circ + \eps]; \Xi^{(\udelta)}\right)$, where
\begin{equation}\label{E:radius}
    \eps \le \min \set{\frac{\Lambda_{(\Bx_\circ, \Bv_\circ)}}{M}, \frac{1}{K_{\Xi^{(\udelta)}}}},
\end{equation}
uniquely exists. Here, 
\begin{equation*} M \Def \sup \set{\norm{\mathcal B(t, \Bx, \Bv;u)}_{\ell^2(\RR^4)}: (t, \Bx, \Bv) \in [t_\circ, T] \times \Xi^{(\udelta)}},\end{equation*}
$K_{\Xi^{(\udelta)}}$ is local Lipschitz constant of function $\mathcal B$ and 
\small{\begin{equation}\label{E:distance}
    \Lambda_{(\Bx_\circ, \Bv_\circ)} \Def \inf \set{ \norm{(\Bx_\circ, \Bv_\circ)- (\By, \Bw)}_{\ell^2}:  (\By,\Bw) \in \partial \Xi^{(\udelta)}},
\end{equation}}
\normalsize
i.e., the distance of the initial point to the boundaries of set $\Xi^{(\udelta)}$. It noteworthy that by compactness of $\Xi^{(\udelta)}$, the infimum term in \eqref{E:distance} is well-defined. This completes the proof of the well-posedness of a local solution. To prove Theorem \ref{T:main}, we need to construct the solution beyond the compact set $\Xi^{(\udelta)}$. This is the subject of the next theorem.
\begin{theorem}[\textbf{Extension of solution}]\label{T:extension} The solution to the system \eqref{E:dynamics} exists and is unique over the domain $\domain$ (defined as in \eqref{E:state_space}).
\end{theorem}
\begin{proof}[\textbf{of Theorem \ref{T:extension}}] Due to the structure of the problem and the explicit form of $\Xi^{(\udelta)}$ (see \eqref{E:compactSet}), the extension of the solution to the entire domain is non-trivial. Canonically, we start with investigating the extension of the solution to the boundaries $\partial \Xi^{(\udelta)}$ of $\Xi^{(\udelta)}$. 
To do so, we define 
\begin{equation*}
    \tau^{(\udelta)} \Def \inf \set{t > t_\circ: (\Bx^{(\udelta)}(t), \Bv^{(\udelta)}(t)) \in \partial\Xi^{(\udelta)}},
\end{equation*}
as the first time that the solution approaches the boundaries. Therefore, as $t \to \tau^{(\udelta)}$, $\eps \to 0$ (see \eqref{E:radius}).
This implies that, on the interval $\mathcal I \Def \left[t_\circ, \min \set{\tau^{(\udelta)}, T}\right)$, i.e., before approaching the boundaries of $\Xi^{(\udelta)}$, the previous discussion of existence and uniqueness of the solution is directly applicable. 
If $\tau^{(\udelta)} > T$, then the existence of a unique solution on $[0, T] \times \mathcal S$ follows immediately since the solution remains in $int(\Xi^{(\udelta)})$ for all $t \in [t_\circ, T]$.\\
So, it suffices to consider $\tau^{(\udelta)} < T$. It should be observed that, in this case, since the solution $(\Bx^{(\udelta)}, \Bv^{(\udelta)})$ is uniformly continuous over the interval $\mathcal I$, it can be extended continuously to $t = \tau^{(\udelta)}$, and hence the solution over $[t_\circ, \tau^{(\udelta)}]$ is well-defined. 
To extend the solution $(\Bx^{(\udelta)},\Bv^{(\udelta)})$, its behavior at the boundaries, i.e., at $t = \tau^{(\udelta)}$ needs to be understood carefully.
Let $0 < \odelta < \udelta$ and by the definition \eqref{E:compactSet},  $\Xi^{(\udelta)} \subset \Xi^{(\odelta)}$ and hence, 
\begin{equation}\label{E:distance_compact}
    \Bx^{(\udelta)}(\tau^{(\udelta)}) \in int (\pi_{\Bx}(\Xi^{(\odelta)})),
\end{equation}
where $\pi_\Bx(\Xi^{(\odelta)})$ is the projection of $\Xi^{(\odelta)}$ on the $\Bx$-coordinates. On the other hand, if the solution $v^{(\udelta)}$ hits the boundary $\set{0, \bar v}$ then beyond the boundary the solution is physically infeasible as either the velocity is negative or will surpass the maximum possible. Mathematically, this implies that if $v^{(\udelta)}(\tau^{(\udelta)}) \in \partial \pi_\Bv(\Xi^{(\udelta)})$, then $\eps = 0$ (see \eqref{E:radius}) and the solution cannot be extended in the same way. The next result shows the velocity function remains in the admissible range and hence we can extend the solution beyond the boundaries.
\begin{lemma}[\textbf{Maximum Principle}]\label{T:boundary}
Fix $\udelta >0$. Define
\begin{equation}
    \mathcal N \Def \set{(\Bx, \Bv) \in [\udelta, \nicefrac{1}{\udelta}] \times (0, \bar v)}.
\end{equation}
Then, $\mathcal N$ is invariant over $[t_\circ, \tau^{(\udelta)}]$ along the trajectory $(\Bx^{(\udelta)}, \Bv^{(\udelta)})$ starting from $(\Bx_\circ, \Bv_\circ)$. In other words, starting from $\mathcal N$ the solution remains in this set for $t \le \tau^{(\udelta)}$.
\end{lemma}
\begin{proof}[\textbf{of Lemma \ref{T:boundary}}]. The proof is by contradiction. In particular, let us suppose on the contrary that the velocity function vanishes, i.e.,
\begin{equation}\label{E:boundary_zero}
    \lim_{t \to \tau^{(\udelta)}} v^{(\udelta)}(t) = 0.
\end{equation}
In this case, the continuity of the solutions (see the Definition \ref{def:solution}), the definition of $\tau^{(\udelta)}$ and the fact that $v^{(\udelta)} >0$ on $[t_\circ, \tau^{(\udelta)})$ imply that there exists a $\check \tau < \tau^{(\udelta)}$, such that $\dot v^{(\udelta)} <0$ on $[\check \tau,  \tau^{(\udelta)})$ (i.e., $v^{(\udelta)}$ decreases on this interval). Considering the dynamics in \eqref{E:accel} and by the definition of the control, $u - v^{(\udelta)} >0$ for sufficiently small $v^{(\udelta)}$. On the other hand, by \eqref{E:boundary_zero}, and since $x^{(\udelta)}_\ell - x^{(\udelta)} \ge \udelta$ on $\mathcal N$, we may choose $\check \tau$ properly such that $v^{(\udelta_)} < \nicefrac{(x^{(\udelta)}_\ell - x^{(\udelta)})}{\tau_s}$. 
Hence, in both cases, for sufficiently small $v^{(\udelta)}$, $\dot v^{(\udelta)} >0$ which contradicts \eqref{E:boundary_zero}. \\
\noindent The case of $v^{(\udelta)}(\tau^{(\udelta)}) = \bar v$ is precluded by the range of control function $u$.
This concludes the proof of Lemma \ref{T:boundary}\end{proof}
\noindent Lemma \ref{T:boundary} and \eqref{E:distance_compact} in particular ensures that 
\begin{equation*}\label{E:veloc_compact}
   (\Bx^{(\udelta)}(\tau^{(\udelta)}),\Bv^{(\udelta)}(\tau^{(\udelta)})) \in int \Xi^{(\odelta)},
\end{equation*}
where $\odelta < \udelta$ and hence can be considered as the initial value for the dynamic model \eqref{E:dynamics} in the larger set $\Xi^{(\odelta)}$. Therefore, by \eqref{E:radius} the solution can be extended uniquely over this compact set. In particular, we should be able to continue the solution in the same way.
 Formally, to complete the proof of Theorem \ref{T:extension}, we note that 
\begin{equation*} \label{E:cover} \domain = \bigcup \set{\Xi^{(\udelta)}: \udelta>0},
\end{equation*} 
i.e., the increasing collection $\set{\Xi^{(\udelta)}: \udelta >0}$ covers the state-space $\mathcal S$. Collecting all together, the solution of dynamical model \eqref{E:dynamics} can be extended to $\mathcal S$ in a unique way. That is, the solution is well-posed before any collision happens. In other words, formally, we have proven the following result:
\begin{corollary}\label{T:collision_time_cor}
Let us define 
\begin{equation} \label{E:collision_time}
    \mathcal T \Def \inf \set {t >t_\circ: x_\ell(t) -x (t) =0}
\end{equation}
i.e., $\mathcal T$ is the collision time. Then, the dynamics \eqref{E:vec_dynamics} has a unique solution on $[0, \mathcal T)$. 
\end{corollary}
This completes the proof of Theorem \ref{T:extension}.  \end{proof}
\noindent So far, the well-posedness is proven for any time before the collision. From here on since the solution does not depend on any set $\Xi^{(\udelta)}$, we drop $\udelta$ and the solution is denote by $t \mapsto (\Bx(t), \Bv(t))$. To complete the prove of Theorem \ref{T:main}, hence, we need to show that the collision does not happen over $[t_\circ, T]$; in other words, we need to show that $\mathcal T > T$.  \\
To do so, we will use the properties of the dynamic model \eqref{E:dynamics} to prove that the distance cannot be less than an explicit threshold. More precisely, we have that  
{\begin{multline*}
       \dot v_\ell(t) - \dot v(t)  \ge \sa_\ell(t) - k_v \frac{v_\ell(t) - v(t)}{(x_\ell(t) - x(t))^2}\\
       \quad - k_d (x_\ell(t)- x(t) - \tau_s v(t)). 
\end{multline*}}
\normalsize
Integrating both sides, we can write
\begin{equation*}
    \begin{split}
        & v_\ell(t) -  v(t)  \\
         & \quad  \ge (v_{\ell, \circ} - v_\circ) + \int_0^t \sa_\ell(s) ds - k_d \int_0^t (x_\ell(s) - x(s) - \tau_s v(s)) ds \\
        & \quad \quad - k_v \int_0^t \frac{v_\ell - v}{(x_\ell - x)^2} ds\\
       & \quad = (v_{\ell, \circ} - v_\circ) + \int_0^t \sa_\ell(s) ds - k_d \int_0^t (x_\ell(s) - x(s) - \tau_s v(s)) ds \\
        &\quad \quad+ k_v \left(\frac{1}{x_\ell(t) - x(t)} - \frac{1}{x_{\ell, \circ} - x_\circ}\right).
    \end{split}
\end{equation*}
Now, let us define $ h(t) \Def x_\ell(t) - x(t)$. Then, we have 
\begin{equation*}\label{E:acc_difference}
 \begin{split}
        \dot h(t) & \ge v_{\ell,\circ} - v_\circ +\int_0^t \sa_\ell(s) ds  - k_d \int_0^t h(s) ds + k_d \tau_s\int_0^t v(s)ds \\
        & \quad + k_v \left(\frac{1}{h(t)} - \frac{1}{h_\circ} \right).
    \end{split}
\end{equation*}

\noindent Now, we continue the proof with contradiction. In particular, we assume that $\mathcal T \le T$, i.e., the collision happens before the time horizon $T$. Then, there exists a time $\check t$ and $\underline h >0$ (see \eqref{E:h_lower_bar}) defined by
\begin{equation}\label{E:tcheck}
    \check t \Def \sup \set{t \in (t_\circ, \mathcal T): x_\ell - x > \underline h}.
\end{equation}
In other words, $x_\ell - x \le \underline h$ on $[\check t, \mathcal T)$. 
By Corollary \ref{T:collision_time_cor} the solution $t \mapsto (\Bx(t), \Bv(t))$ exists over $[0, \mathcal T)$. Noting that $h:[0, \mathcal{T}) \to \mathbb{R}_+$ is continuous by the Definition \ref{def:solution}, it is bounded over the compact set $[0, \check t]$. In addition, by the definition of $\check t$, $h(t) \le \underline h$ on $[\check t, \mathcal T)$. Putting all together, function $h$ is integrable over the time interval  $[0, \mathcal{T})$. In other words, there exists $H>0$ such that 
\begin{equation*} 
\int_{t_\circ}^{\mathcal T} h(s) ds \leq H.
\end{equation*}
In addition, $ v_{\ell}(t)=v_{\ell,\circ}+ \int_{t_\circ}^t \sa_{\ell}(s) ds  \ge 0$ and $k_d \tau_s \int_0^t v(s)ds \ge 0$, $\forall t \in [0, \mathcal{T})$. Thus, we have,  

\begin{equation}\label{E:doth}
        \dot{h}(t) \geq - v_\circ - Hk_d -\frac{k_v}{h_\circ} + \frac{k_v}{h(t)}.
\end{equation}
Let us define 
\begin{align}
\underline h \Def \min \set{\frac{k_v}{v_\circ+Hk_d+\frac{k_v}{h_{\circ}}}, h_\circ}=\frac{k_v}{v_\circ+Hk_d+\frac{k_v}{h_{\circ}}},
\label{E:h_lower_bar}
\end{align}
where $h_\circ \Def x_{\ell, \circ} - x_\circ$. Then, \eqref{E:h_lower_bar} and \eqref{E:doth} imply that 
\begin{equation*}
    \dot h(t) \ge 0, \forall t \in [\check t, \mathcal T).
\end{equation*}
In other words, 
\begin{equation*}
    \label{E:non-collision}
    \inf_{t \in [\check t, \mathcal T)} h(t) \ge \underline h,
\end{equation*}
which by continuity of $h(\cdot)$ on $[t_\circ, \mathcal T)$ is contradiction to the assumption $\mathcal T \le T$. Thus, $\mathcal T >T$ and by Corollary \ref{T:collision_time_cor}, the solution of \eqref{E:dynamics} uniquely exists on $[t_\circ, T]$. 

Finally, if the solutions $(x_\ell, x, v_\ell, v)$ exist in the sense of the Definition \ref{def:solution}, then they belong to a class of absolutely continuous functions and hence continuous and bounded on $[0, T]$. In addition, $\dot v(t) = \mathcal A(\Bx, \Bv; u)$ a.e. and $\norm{\mathcal A(\Bx, \Bv; u)}_{L^\infty(0, T)} \le \norm{u - v}_{L^\infty(0, T)} < \infty$. This in particular implies that $(\Bx, \Bv) \in W^{1, \infty}(0, T; \mathcal S)$. 

\begin{figure}[ht] 
  \begin{subfigure}[b]{0.5\linewidth}
    \centering
    \includegraphics[width=0.85\linewidth]{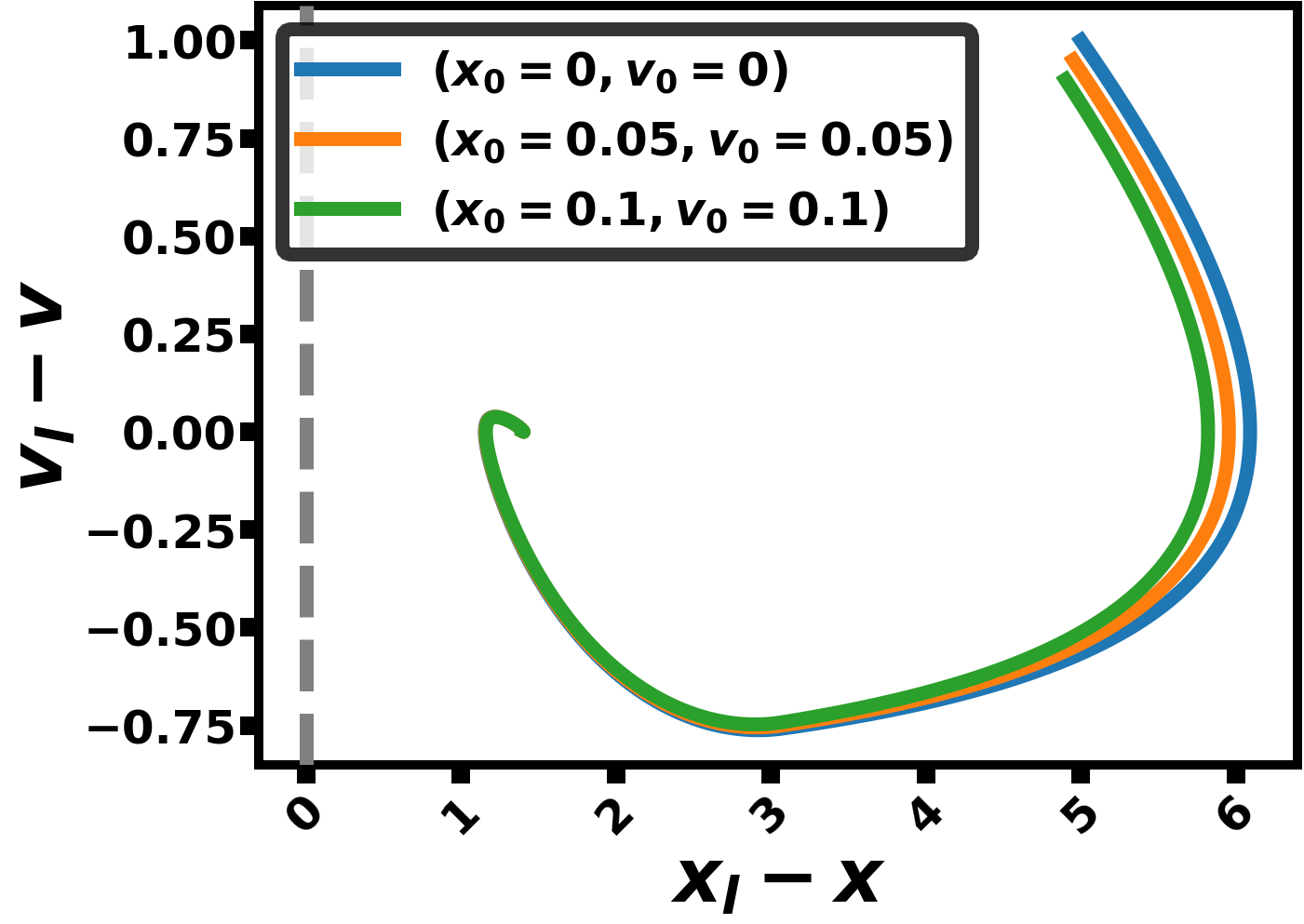} 
    \caption{\footnotesize{Fixed $x_{\ell, \circ} = 5, v_{\ell, \circ} = 1$}}
    \label{fig: initial_condition_close_case1} 
  \end{subfigure}
  \begin{subfigure}[b]{0.5\linewidth}
    \centering
    \includegraphics[width=0.85\linewidth]{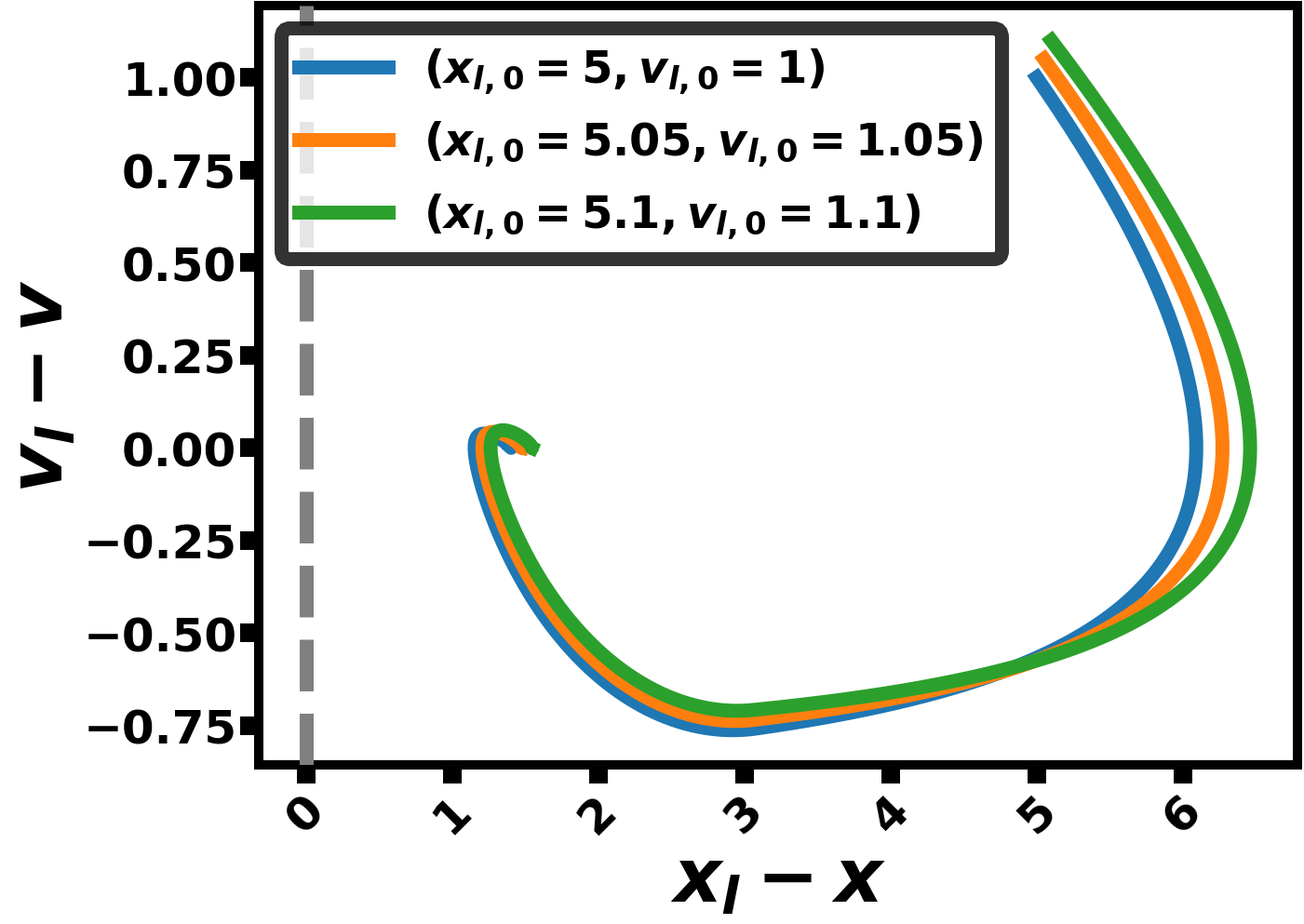} 
    \caption{\footnotesize{Fixed $x_\circ = 0, v_\circ = 0$.}}
    \label{fig: initial_condition_close_case2} 
  \end{subfigure} 
  \vspace{-10pt}
  \caption{\scriptsize{Different cases with initial conditions being close (0, 0.05, and 0.1). Parameters same as Fig \ref{fig: proposed}.}}
  \label{fig:stability} 
  \vspace{-10pt}
\end{figure}

The last step is to investigate the continuous dependence of the solution on initial data in uniform topology. By construction of the solution, in particular boundedness of $v_\ell$ and $v$, the boundedness of $x_\ell$ and $x$ over $[0, T]$ and subsequently the Lipschitz continuity of $(\Bx, \Bv) \mapsto \mathcal B(t, \Bx, \Bv;u)$ follows. Let $\varphi \Def (\Bx, \Bv)$ and $\psi \Def (\By, \Bw)$ be the solutions of the dynamical model \eqref{E:dynamics} with respect to the initial data $\phi_\circ = (\Bx_\circ, \Bv_\circ)$ and $\psi_\circ = (\By_\circ, \Bw_\circ)$, respectively. Then, using the Lipschitz continuity and Gronwall's inequality, we have 
\begin{equation}
    \norm{\varphi - \psi}_{L^\infty(0, T)} \le \norm{\varphi_\circ - \psi_\circ}_{\ell^2(\RR^4)} \exp\lb C_L T \rb,
\end{equation}
where $C_L$ is Lipschitz constant of function $\mathcal B$. 
This completes the proof of Theorem \ref{T:main} (Fig. \ref{fig:stability} illustrates the concept in two cases).
\end{proof}
\begin{remark}[\textbf{Collision and Generalization}] \label{R:generalization}
Proof of Theorem \ref{T:main}, particularly \eqref{E:non_collision}, shows that in the dynamics \eqref{E:dynamics}, the collision does not happen. Since the velocity $v_\ell$ and acceleration $\sa_\ell$ are considered as general time-dependent dynamics, the result can be immediately extended to the case of dynamical model \eqref{E:general}. \hfill $\diamond$
\end{remark}
\subsection{\textbf{Perturbation and Stability}}\label{S:perturbation}
In this section, we study the behavior of the dynamical model \eqref{E:dynamics} in the presence of small perturbation in the leader's velocity. The main goal of such studies is two-fold. Small perturbation method provides a framework to study the stable behavior of the solution along the trajectory. In addition, no physical system is isolated from noise. Therefore, as the behavior of the perturbed systems is unknown, proving that the perturbed system remains close to the original dynamic for small perturbations, provides a proper estimation on the behavior of the perturbed dynamics.
In particular, we consider
\begin{equation}\label{E:perturbed_dynamics}
    \dot v_\ell(t) = \sa_\ell(t) + \eps g(t),
\end{equation}
for some $\eps >0$, where $g(\cdot)$ is a perturbation function, and the solution of the perturbed dynamics\eqref{E:perturbed_dynamics} is denoted by $(x^\eps_\ell, x^\eps, v^\eps_\ell, v^\eps)$; cf. \eqref{E:dynamics}.
\begin{theorem}
    Let $g \in C([0, T];\RR)$ be the perturbation function which satisfies
    \begin{equation} \label{E:bound_guarantee}
    \norm{g}_{L^1(0, T)} \le \tfrac {1}{\eps_\circ} \left(\bar v - v_{\ell, \circ} - \norm{\sa_\ell}_{L^1(0, T)}\right) ,
    \end{equation} 
    for some $\eps_\circ >0$. Let $\xi^\eps \Def x^\eps_\ell - x^\eps$ and $\zeta^\eps \Def v^\eps_\ell - v^\eps$, where $(x_\ell^\eps, x^\eps, v_\ell^\eps, v^\eps)$ is the solution to the perturbed dynamics \eqref{E:perturbed_dynamics} for any $\eps \le \eps_\circ$. Similarly, let $\xi \Def x_\ell- x$ and $\zeta \Def v_\ell- v$ where $(x_\ell, x, v_\ell,v)$ is the solution to the dynamical model \eqref{E:dynamics}. Then,
    \begin{equation}\label{E:converg}
        \norm{(\xi^\eps, \zeta^\eps)^\sT - (\xi, \zeta)^\sT}_{L^\infty(0, T)} \to 0 \quad, \text{as $\eps \to 0$}.
    \end{equation}
\end{theorem}
\begin{proof} 
Eq. \eqref{E:bound_guarantee} ensures that the leading velocity in the presence of the perturbation does not violate the maximum speed $\bar v$. In particular, by the maximum principle in Lemma \ref{T:boundary}, we conclude uniform boundedness of the form 
    \begin{equation}\label{E:uniform}
        \sup_{\eps >0}\sup_{t \in [0,T]} \abs{\zeta^\eps (t)} < \bar v.
    \end{equation}
    Furthermore, using the dynamics of \eqref{E:perturbed_dynamics}, we have that 
    \small{
    \begin{multline*}
         \zeta^\eps(t) = \zeta_\circ - \int_0^t \mathcal A(\xi^\eps(s), \zeta^\eps(s); u) ds + \int_0^t \sa_\ell(s) ds + \eps \int_0^t g(s) ds.
    \end{multline*}
    }
\normalsize
    By \eqref{E:bound_guarantee} and consequently dominated convergence theorem, the last integral term vanishes as $\eps \to 0$. Hence, by Theorem \ref{T:main}, we can pass the limit as $\eps \to 0$ and we get
   \small{{\begin{equation*}
        \sup_{t \in [0, T]} \Big|\zeta^\eps(t) - 
         \left(\zeta_\circ + \int_0^t \sa_\ell(s) ds - \int_0^t \mathcal A(\xi^\eps(s), \zeta^\eps(s); u) ds \right) \Big| \to 0 .
    \end{equation*}
    }}
    \normalsize
    A similar argument provides a similar result for $\xi^\eps$.
    Therefore, to conclude that $(\xi^\eps, \zeta^\eps) \to (\xi, \zeta)$ in uniform topology, we need to show that the limit of $\zeta^\eps$ and $\xi^\eps$ as $\eps \to 0$ exists and then by continuity of $\mathcal A$ the result follows. To this end, let us define $\mathcal M \Def \set{(\xi^\eps,\zeta^\eps): \eps>0} \subset C([0, T]; \RR^2)$ such that $(\xi^\eps, \zeta^\eps)$ is the solution of the dynamical model \eqref{E:perturbed_dynamics} (when written in the difference form). In addition, \eqref{E:uniform} invokes the Arzela-Ascoli theorem, and hence $\mathcal M$ is totally bounded (precompact) in the Banach space $(C([0, T]; \RR^2), \mathcal T_{\norm{\cdot}_\infty})$; the uniform topology. Therefore, any sequence $\set{(\xi^{\eps_n}, \zeta^{\eps_n}): n \in \NN}$ in which $\eps_n \to 0$ as $n \to \infty$, has a convergent subsequence with the limiting function $(\hat \xi, \hat \zeta) \in C([0, T]; \RR^2)$ which clearly satisfies the dynamical model \eqref{E:dynamics}. The existence of a unique solution by Theorem \ref{T:main}, on the other hand, proves the claimed result \eqref{E:converg}.
\end{proof}

\begin{figure}[h!]
    \centering
    \includegraphics[width=0.9\linewidth]{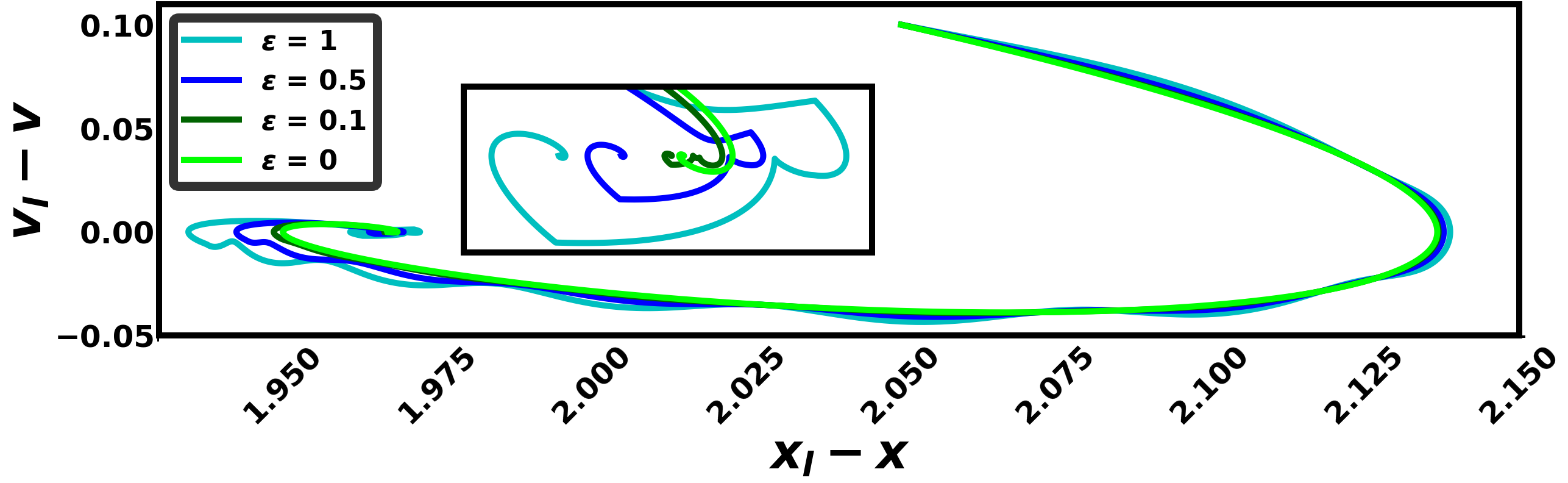}
    \caption{\scriptsize{Perturbation of Leading Vehicle's Velocity. Parameters same as Fig. \ref{fig: proposed}, except $T = 60$. Here, $(x_{\ell, \circ}, x_\circ, v_{\ell, \circ}, v_\circ) = (2.05,0,1.4,1.3)$. The subplot focuses on near the sink: $\sup_{t \in [0, T]} |\zeta^\eps(t) - \zeta(t)| = 49 \times 10^{-4}$ when $\eps = 1$. 
    }}
    \label{fig:perturbation}
\end{figure}
Figure \ref{fig:perturbation} illustrates the perturbation of the leading vehicle with constant velocity with a perturbation function constructed from various functions in increasing, decreasing, and periodic form to capture the behavior of the perturbed dynamics for several $\eps$ values. 
\subsection{\textbf{Primary Estimations of Solution}}\label{S:Estimates}
\begin{figure}[h!]
    \centering
    \includegraphics[width=0.9\linewidth]{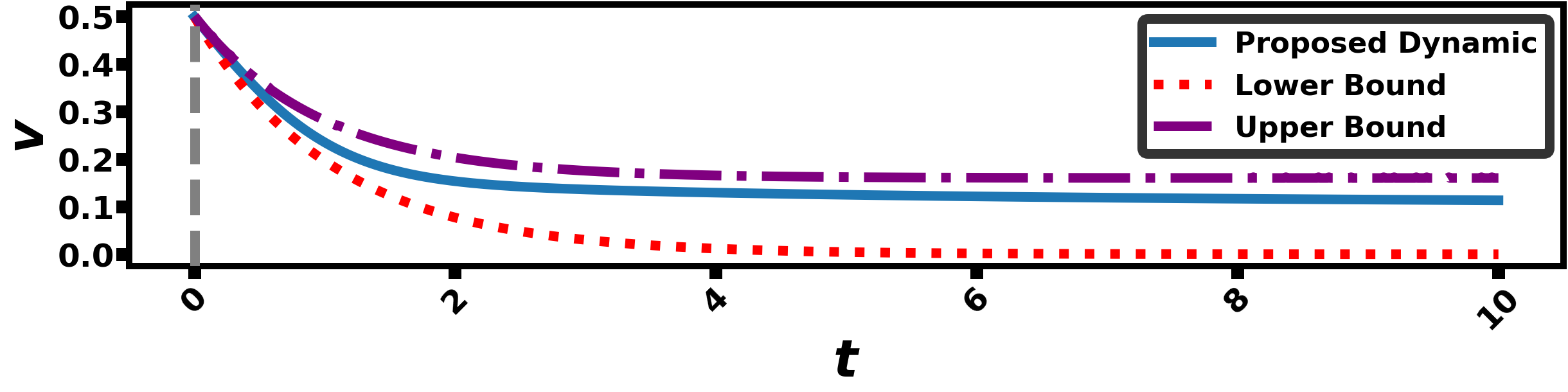}
    \vspace{-10pt}
    \caption{\scriptsize{Bounds on {follower} vehicle's velocity. Parameter same as Fig. \ref{fig: proposed} except $T = 10$. Initial Conditions: $(x_{\ell, \circ}, x_\circ, v_{\ell, \circ}, v_\circ) = (1,0,0.1,0.5)$.}} 
    \label{fig:bound_v} 
    \vspace{-10pt}
\end{figure}

In this section, we derive some estimates on the solutions of the dynamics \eqref{E:accel} along the trajectories; i.e., time-dependent bounds which will provide more information about the behavior of the solution.
Consider the dynamics governed by \eqref{E:dynamics}, the following vehicle's acceleration is bounded below by 
    \begin{equation}\label{E:initial_v_upperbound}
    \begin{split}
   \dot{v} & = \min\set{k_v \frac{( v_\ell -  v)}{(x_\ell- x)^2} + k_d (x_\ell- x - \tau_s v), k (u - v)}\\
    & \ge -\max\set{\tfrac{k_v}{\underline h ^2}+k_d \tau_s, k} v,
\end{split}\end{equation}
where $\underline h$ as in \eqref{E:h_lower_bar} and control term $u \in C((0, T); [\underline u, \bar u])$.
Considering \eqref{E:initial_v_upperbound}, one can prove by contradiction that for any $t \in [0, T]$ and the initial velocity $v_\circ$ 
\begin{align}
\label{E: velocity_lower_bound}
    v(t)\geq \underline {\mathcal V} (t) \Def v_\circ e^{-\max\set{\tfrac{k_v}{\underline h ^2}+{k_d\tau_s}, k}t}.
\end{align}
Next, we derive the corresponding upper bound.
Using Lemma \ref{T:boundary} and \eqref{E: velocity_lower_bound} and some algebraic manipulations, for any $t \in [0,T]$, we have that
\small{{\begin{equation}
\label{E: h_upper_bound}
    h(t) \leq \bar{h}(t) \Def  h_\circ + \bar{v}t +v_\circ  
\nicefrac{\lb e^{-\max\set{\tfrac{k_v}{\underline h^2}+{k_d\tau_s}, k}t} -1 \rb}{\max\set{\tfrac{k_v}{\underline h^2}+{k_d\tau_s}, k}}.
\end{equation}
}}
\normalsize
Using \eqref{E:accel}, Lemma \ref{T:boundary}, \eqref{E:non_collision} and \eqref{E: h_upper_bound}, we can write
\begin{multline*}
    \dot{v}(t) \le \min \set{-k_d \tau_s v(t) + k_d \bar{h}(t) + k_v \frac{\bar{v}}{\underline h^2}, -kv(t) + ku(t)},
\end{multline*}
and hence, for any $t \in [0, T]$
\begin{multline}\label{E:v_upperbound}
    v(t) \leq  \bar{\mathcal V}(t) \Def \min\Bigg\{k \int_0^t e^{k(s-t)} u(s) \, ds  + v_{\circ} e^{-kt}, \\ 
     \int_0^t e^{k_d \tau_s (s-t)} \Big( k_d \bar{h}(s) + k_v \frac{\bar{v}}{\underline h^2} \Big)\, ds  + v_{\circ} e^{-k_d \tau_s t}\Bigg\}. 
\end{multline}
Figure \ref{fig:bound_v} illustrates the {bounds} \eqref{E: velocity_lower_bound} and \eqref{E:v_upperbound}. We have used a simple calibration to tighten the bounds. The optimal calibration of the constants {determines the tightness of these bounds and} is outside the scope of this paper. 
\section{Conclusions and {future work}}
\vspace{-4pt}
A nonlinear dynamical model was considered in this paper. We rigorously proved the well-posedness in some functional spaces. The bounds along the trajectories of the solutions are established. It is evident that depending on the initial conditions, calibration of parameters can remarkably improve these bounds. A comprehensive study of such calibration will be a future work. Finally, we have shown that under a small perturbation, the behavior of the dynamical model remains close to the original dynamics. This shows the stability of the solution as well as a proper estimation of the solution in the presence of a small perturbation.
\vspace{-4pt}

\bibliographystyle{IEEEtran}
\bibliography{reference}

\begin{thebibliography}{10}
\providecommand{\url}[1]{#1}
\csname url@samestyle\endcsname
\providecommand{\newblock}{\relax}
\providecommand{\bibinfo}[2]{#2}
\providecommand{\BIBentrySTDinterwordspacing}{\spaceskip=0pt\relax}
\providecommand{\BIBentryALTinterwordstretchfactor}{4}
\providecommand{\BIBentryALTinterwordspacing}{\spaceskip=\fontdimen2\font plus
\BIBentryALTinterwordstretchfactor\fontdimen3\font minus \fontdimen4\font\relax}
\providecommand{\BIBforeignlanguage}[2]{{%
\expandafter\ifx\csname l@#1\endcsname\relax
\typeout{** WARNING: IEEEtran.bst: No hyphenation pattern has been}%
\typeout{** loaded for the language `#1'. Using the pattern for}%
\typeout{** the default language instead.}%
\else
\language=\csname l@#1\endcsname
\fi
#2}}
\providecommand{\BIBdecl}{\relax}
\BIBdecl

\bibitem{treiber2013traffic}
M.~Treiber and A.~Kesting, ``Traffic flow dynamics,'' \emph{Traffic Flow Dynamics: Data, Models and Simulation, Springer-Verlag Berlin Heidelberg}, 2013.

\bibitem{pipes1953operational}
L.~A. Pipes, ``An operational analysis of traffic dynamics,'' \emph{Journal of applied physics}, vol.~24, no.~3, pp. 274--281, 1953.

\bibitem{newell1961nonlinear}
G.~F. Newell, ``Nonlinear effects in the dynamics of car following,'' \emph{Operations research}, vol.~9, no.~2, pp. 209--229, 1961.

\bibitem{gazis1961nonlinear}
D.~C. Gazis, R.~Herman, and R.~W. Rothery, ``Nonlinear follow-the-leader models of traffic flow,'' \emph{Operations research}, vol.~9, no.~4, pp. 545--567, 1961.

\bibitem{treiber2017intelligent}
M.~Treiber and A.~Kesting, ``The intelligent driver model with stochasticity-new insights into traffic flow oscillations,'' \emph{Transportation Research Procedia}, vol.~23, pp. 174--187, 2017.

\bibitem{bando1998analysis}
M.~Bando, K.~Hasebe, K.~Nakanishi, and A.~Nakayama, ``Analysis of optimal velocity model with explicit delay,'' \emph{Physical Review E}, vol.~58, no.~5, p. 5429, 1998.

\bibitem{wang2023car}
Z.~Wang, Y.~Shi, W.~Tong, Z.~Gu, and Q.~Cheng, ``Car-following models for human-driven vehicles and autonomous vehicles: A systematic review,'' \emph{Journal of transportation engineering, Part A: Systems}, vol. 149, no.~8, p. 04023075, 2023.

\bibitem{zha2023survey}
Y.~Zha, J.~Deng, Y.~Qiu, K.~Zhang, and Y.~Wang, ``A survey of intelligent driving vehicle trajectory tracking based on vehicle dynamics,'' \emph{SAE International journal of vehicle dynamics, stability, and NVH}, vol.~7, no. 10-07-02-0014, 2023.

\bibitem{milanes2013cooperative}
V.~Milan{\'e}s, S.~E. Shladover, J.~Spring, C.~Nowakowski, H.~Kawazoe, and M.~Nakamura, ``Cooperative adaptive cruise control in real traffic situations,'' \emph{IEEE Transactions on intelligent transportation systems}, vol.~15, no.~1, pp. 296--305, 2013.

\bibitem{kerner2015failure}
B.~S. Kerner, ``Failure of classical traffic flow theories: a critical review,'' \emph{e \& i Elektrotechnik und Informationstechnik}, vol.~7, no. 132, pp. 417--433, 2015.

\bibitem{jafaripournimchahi2022integrated}
A.~Jafaripournimchahi, Y.~Cai, H.~Wang, L.~Sun, J.~Weng \emph{et~al.}, ``Integrated-hybrid framework for connected and autonomous vehicles microscopic traffic flow modelling,'' \emph{Journal of Advanced Transportation}, vol. 2022, 2022.

\bibitem{chandler1958traffic}
R.~E. Chandler, R.~Herman, and E.~W. Montroll, ``Traffic dynamics: studies in car following,'' \emph{Operations research}, vol.~6, no.~2, pp. 165--184, 1958.

\bibitem{nick2022near}
H.~Nick Zinat~Matin and R.~B. Sowers, ``Near-collision dynamics in a noisy car-following model,'' \emph{SIAM Journal on Applied Mathematics}, vol.~82, no.~6, pp. 2080--2110, 2022.

\bibitem{dellemonache2019pardalos}
M.~L. Delle~Monache, T.~Liard \emph{et~al.}, ``Feedback control algorithms for the dissipation of traffic waves with autonomous vehicles,'' \emph{Computational Intelligence and Optimization Methods for Control Engineering}, pp. 275--299, 2019.

\bibitem{gong2023well}
X.~Gong and A.~Keimer, ``On the well-posedness of the “bando-follow the leader” car following model and a time-delayed version,'' \emph{Networks and Heterogeneous Media}, vol.~18, no.~2, pp. 775--798, 2023.

\bibitem{matin2023existence}
H.~N.~Z. Matin and M.~L.~D. Monache, ``On the existence of solution of conservation law with moving bottleneck and discontinuity in flux,'' \emph{arXiv preprint arXiv:2310.00537}, 2023.

\bibitem{van2006impact}
B.~Van~Arem, C.~J. Van~Driel, and R.~Visser, ``The impact of cooperative adaptive cruise control on traffic-flow characteristics,'' \emph{IEEE Transactions on intelligent transportation systems}, vol.~7, no.~4, pp. 429--436, 2006.

\bibitem{milanes2014modeling}
V.~Milan{\'e}s and S.~E. Shladover, ``Modeling cooperative and autonomous adaptive cruise control dynamic responses using experimental data,'' \emph{Transportation Research Part C: Emerging Technologies}, vol.~48, pp. 285--300, 2014.

\bibitem{shladover2012impacts}
S.~E. Shladover, D.~Su, and X.-Y. Lu, ``Impacts of cooperative adaptive cruise control on freeway traffic flow,'' \emph{Transportation Research Record}, vol. 2324, no.~1, pp. 63--70, 2012.

\bibitem{matin2023near}
H.~N.~Z. Matin and M.~L. Delle~Monache, ``Near collision and controllability analysis of nonlinear optimal velocity follow-the-leader dynamical model in traffic flow,'' in \emph{62nd IEEE Conference on Decision and Control (CDC)}.\hskip 1em plus 0.5em minus 0.4em\relax IEEE, 2023, pp. 8057--8062.

\bibitem{wang2024hierarchical}
H.~Wang, Z.~Fu, J.~Lee, H.~N.~Z. Matin \emph{et~al.}, ``Hierarchical speed planner for automated vehicles: A framework for lagrangian variable speed limit in mixed autonomy traffic,'' \emph{arXiv preprint arXiv:2402.16993}, 2024.

\bibitem{davis2003modifications}
L.~Davis, ``Modifications of the optimal velocity traffic model to include delay due to driver reaction time,'' \emph{Physica A: Statistical Mechanics and its Applications}, vol. 319, pp. 557--567, 2003.

\bibitem{stern2018dissipation}
R.~E. Stern, S.~Cui, M.~L. Delle~Monache \emph{et~al.}, ``Dissipation of stop-and-go waves via control of autonomous vehicles: Field experiments,'' \emph{Transportation Research Part C: Emerging Technologies}, vol.~89, pp. 205--221, 2018.

\bibitem{adams2003sobolev}
R.~A. Adams and J.~J. Fournier, \emph{Sobolev spaces}.\hskip 1em plus 0.5em minus 0.4em\relax Elsevier, 2003.

\end{thebibliography}
\end{document}